\documentclass[10pt]{dcds2}

\def\ppages{1--10}

\theoremstyle{plain}
\newtheorem{theorem}{Theorem}
\newtheorem{corollary}{Corollary}

\newtheorem{proposition}{Proposition}

\theoremstyle{definition}
\newtheorem{definition}{Definition}
\newtheorem{remark}{Remark}

\newcommand{\real}{\mathbb{R}}

\newcommand{\Lc}{{\cal L}}
\newcommand{\ov}{\overline}
\newcommand{\D}{{\mathcal{D}}}
\newcommand{\V}{{\mathcal{V}}}
\newcommand{\A}{{\mathcal{A}}}

\def\g{\mathfrak{g}}

\newcommand{\DEL}{\operatorname{DEL}}
\newcommand{\Ad}{\operatorname{Ad}}
\newcommand{\Id}{\operatorname{Id}}

\title[Energy conserving nonholonomic integrators]{Energy conserving nonholonomic integrators}

\thanks{Former address: Instituto de Matem\'aticas y F{\'\i}sica Fundamental,
  Consejo Superior de Investigaciones Cient{\'\i}ficas, Serrano 123, 28006
  Madrid, Spain}

 \email{jcortes@uiuc.edu}
 \urladdr{motion.csl.uiuc.edu/\~{}jorge}


\author[Jorge Cort\'es]{}

\def\cal#1{{\fam2#1}}

\begin{document}
\maketitle

\centerline{\scshape Jorge Cort\'es}
\smallskip

{\footnotesize \centerline{Coordinated Science Laboratory}
  \centerline{University of Illinois at Urbana-Champaign}
  \centerline{1308 W. Main St., IL 61801} \centerline{United States} }

\begin{quote}{\normalfont\fontsize{8}{10}\selectfont
    {\bfseries Abstract.}  We address the problem of constructing
    numerical integrators for nonholonomic Lagrangian systems that
    enjoy appropriate discrete versions of the geometric properties of
    the continuous flow, including the preservation of energy.
    Building on previous work on time-dependent discrete mechanics,
    our approach is based on a discrete version of the
    Lagrange-d'Alembert principle for nonautonomous systems.
\par}
\end{quote}

\section{Introduction}

\renewcommand{\descriptionlabel}[1]
{\hspace{\labelsep}\textsf{#1}}

In the last years Geometric Integration has grown to be a very large
and active area of research, with a rich variety of approaches taken
and topics covered~\cite{BuIs,Sa}.  Among the various viewpoints, the
variational integrators approach has revealed to be very
powerful~\cite{Ma2}.  This point of view is not confined to Lagrangian
and Hamiltonian (conservative) systems, but also admit extensions to
multisymplectic geometry and PDEs, as well as to systems subject to
external forces and dissipation (see~\cite{MaWe} for a recent overview
on the subject).

The treatment of problems with constraints has also been an important issue in
the area.  Holonomic constraints have received a great deal of
attention~\cite{Go,Ja,LeRe,LeSk,Re,Se}, motivated by their presence in
applications such as molecular dynamics and planetary motions.  The treatment
of nonholonomic constraints has also been in the agenda of the Geometric
Integration community (see, for instance,~\cite{McSc,WeMa}).  Following the
variational approach to discrete mechanics, we proposed in~\cite{CoMa} a class
of nonholonomic numerical integrators enjoying discrete versions of some of
the geometric properties of the continuous flow.  These include the evolution
of the symplectic form along the flow, and the fulfillment of a discrete
version of the nonholonomic momentum equation~\cite{BlKrMaMu}, which in the
case of horizontal symmetries gives rise to conservation laws.  However, these
integrators do not preserve the energy, which is a natural conserved quantity
of the continuous flow. This is not surprising, since (fixed time-step)
variational integrators themselves do not preserve the energy either.  A
different approach based on the technique of generating functions is proposed
in~\cite{LeMaSa}.

In this paper, we address the problem of energy conservation building on
previous derivations on time-dependent discrete mechanics and extended
variational integrators~\cite{KaMaOr,LeMa,MaWe}. Our main contribution is the
construction of extended nonholonomic integrators derived from a discrete
version of the Lagrange-d'Alembert principle for nonautonomous systems.  We
focus on investigating the relationships between the discrete and the
continuous mechanics.  The special feature of these integrators is that, in
addition to inheriting good properties with respect to the symplectic form and
the nonholonomic momentum, they also preserve the energy.

The paper is organized as follows.
Section~\ref{se:mechanical-integrators} gives a brief introduction to
mechanical integrators for unconstrained systems and explains the
necessity of allowing variable time steps to design algorithms which
preserve at the same time the symplectic form, the momentum and the
energy. The basic theory on time-dependent variational integrators is
also presented.  In Section~\ref{se:EDLA}, we propose a discrete
version of the Lagrange-d'Alembert principle for nonautonomous
constrained Lagrangian systems.  This principle leads us naturally to
the extended discrete Lagrange-d'Alembert equations, which we term
nonholonomic integrators.  Section~\ref{se:EDLA-geometric-properties}
presents an account of the geometric properties of these integrators,
paying special attention to the energy conservation.  Finally,
Section~\ref{se:conclusions} gives some concluding remarks.

\section{Mechanical integrators}\label{se:mechanical-integrators}

In this section, we briefly introduce some common notions and results
from the literature on Geometric Integration. For further reference,
the reader is referred to~\cite{Ha,Ma,MaPaSh,SaCa}.  Given a
symplectic manifold $(P,\omega)$ and a Hamiltonian function $H: P
\rightarrow \real$, an algorithm $F_h: P \rightarrow P$, $h \in
[0,h_0)$, is called a \emph{symplectic integrator} if each $F_h: P
\rightarrow P$ is a symplectic map; an \emph{energy integrator} if $H
\circ F_h = H$; and a \emph{momentum integrator} if $J \circ F_h = J$,
where $J:P \rightarrow \g^*$ is the momentum map associated with the
action of a Lie group $G$ on $P$.  An algorithm having any of these
properties is called a \emph{mechanical integrator}.

The choice of a specific integrator depends on the concrete problem
under consideration.  For instance, in molecular dynamics simulation,
the preservation of the symplectic form is important for long time
runs, since otherwise one may obtain totally inconsistent solutions.
On the other hand, the exact conservation of momentum first integrals
is essential to problems in attitude control in satellite dynamics,
since this is the basic physical principle driving the reorientation
of the system. However, one is in general prevented from finding
integrators which preserve the three elements at the same time due to
the following result.

\begin{theorem}[\cite{GeMa}]\label{th:Ge-Marsden}
  Consider a Hamiltonian system with a symmetry group $G$ such that
  the dynamics $X_H$ is nonintegrable on the reduced space (in the
  sense that any other conserved quantity is functionally dependent on
  $H$). Assume that a numerical integrator for this system is
  energy-symplectic-momentum preserving and $G$-equivariant. Then, the
  integrator gives the exact solution of the problem up to a time
  reparameterization.
\end{theorem}

Roughly speaking, this result means that obtaining a fixed time step
 energy-symplectic-momentum integrator is the same as
exactly obtaining the continuous flow.  This theoretical obstruction
can be overcome by allowing for varying time steps~\cite{KaMaOr}, as
we will review below.

\subsection{Variational integrators}\label{se:variational-integrators}

Mechanical integrators derived from discrete mechanics have their
origin in the works by Lee, Veselov and others
(see~\cite{Le,MoVe,Ve,Ve2} and references therein).  In the last
years, they have been intensively studied and further developed to
deal with more general situations~\cite{BoSu,KaMaOr,KaMaOrWe,WeMa}.
We briefly review here the main ideas of this approach.  A complete
exposition can be found in the recent overview~\cite{MaWe}.  For the
sake of conciseness, we directly go to the time-dependent case,
without presenting the autonomous situation.

Let $Q$ be an $n$-dimensional manifold, and consider the extended
configuration manifold $\ov{Q} = \real \times Q$.  The extended
discrete Lagrangian state space is $\ov{Q} \times \ov{Q}$, with
canonical projections $\pi_i: \ov{Q} \times \ov{Q} \rightarrow
\ov{Q}$, $i=1,2$.  An \emph{extended discrete path} is a sequence of
points in $\ov{Q}$, i.e. a map $c: \{ 0,\dots, N \} \rightarrow
\ov{Q}$.  We denote $c(k)=(t_k,q_k) \in \ov{Q}$, $k=0,\dots,N$.  Given
a discrete path, the associated discrete curve is $q: \{ t_0,\dots,t_N
\} \rightarrow Q$, $q (t_k) = q_k$.  The extended discrete path space
is defined by
\[
C_d = \{ c: \{ 0,\dots, N \} \rightarrow \ov{Q} \; | \; t_{k+1} > t_k
\, , \; k = 0,\dots,N-1\} \, .
\]
The tangent space $T_cC_d$ to $C_d$ at $c$ is the set of all maps
$\delta c: \{ 0,\dots, N\} \rightarrow T\ov{Q}$ such that
$\tau_{\ov{Q}} \circ \delta c =c$, where $\tau_{\ov{Q}}:T\ov{Q}
\rightarrow \ov{Q}$ denotes the canonical projection.  Consider the
space $(\ov{Q} \times \ov{Q})^2 = \ov{Q} \times \ov{Q} \times \ov{Q}
\times \ov{Q}$ with projections $\sigma_i: (\ov{Q} \times \ov{Q})^2
\rightarrow \ov{Q} \times \ov{Q}$, $i=1,2$.  The extended discrete
second-order manifold of $(\ov{Q} \times \ov{Q})^2$ is defined by $
\ddot{Q}_d = \{ w \in (\ov{Q} \times \ov{Q})^2 \; | \; \pi_2 \circ
\sigma_1 (w) = \pi_1 \circ \sigma_2 (w) \}$. Otherwise said,
$\ddot{Q}_d$ is the set of points $w$ in $(\ov{Q} \times \ov{Q})^2$ of
the form $w = (t_0,q_0,t_1,q_1$, $t_1,q_1,t_2,q_2)$.

An \emph{extended discrete Lagrangian system} is given by a map $L_d:
\ov{Q} \times \ov{Q} \rightarrow \real$.  The \emph{extended action
  sum} $S: C_d \rightarrow \real$ is then defined by,
\begin{equation}\label{actionsum}
  S (c) = \sum_{k=0}^{N-1} L_d (c(k),c(k+1)) =  \sum_{k=0}^{N-1} L_d
  (t_k,q_k,t_{k+1},q_{k+1}) \, . 
\end{equation}

\begin{theorem}[\cite{MaWe}]\label{th:discrete-action}
  Given a $C^k$ extended discrete Lagrangian $L_d:\ov{Q} \times \ov{Q}
  \rightarrow \real$, $k \ge 1$, there exists a unique $C^{k-1}$
  mapping $D_{\DEL} L_d : \ddot{Q}_d \rightarrow T^*\ov{Q}$ and unique
  $C^{k-1}$ one-forms $\Theta_{L_d}^-$ and $\Theta_{L_d}^+$ on $\ov{Q}
  \times \ov{Q}$ such that, for all variations $\delta c \in T_c C_d$
  of $c \in C_d$, \vspace*{-.5cm}
  \begin{multline}
    dS (c) \cdot \delta c = \sum_{k=1}^{N-1} D_{\DEL} L_d
    (t_{k-1},q_{k-1},t_{k},q_{k},t_{k+1},q_{k+1}) \cdot (\delta t_{k},
    \delta q_{k})
    \\
    \qquad + \Theta_{L_d}^+ (t_{N-1},q_{N-1},t_{N},q_{N}) \cdot
    (\delta t_{N-1},\delta q_{N-1},\delta t_{N},\delta q_{N})
    \\
    - \Theta_{L_d}^- (t_{0},q_{0},t_{1},q_{1}) \cdot (\delta
    t_{0},\delta q_{0},\delta t_{1},\delta q_{1}) \, .
  \end{multline}
\end{theorem}

The map $D_{\DEL}L_d$ is called the \emph{extended discrete Euler-Lagrange
  map} and the one-forms $\Theta_{L_d}^+$ and $\Theta_{L_d}^-$ are the
\emph{extended discrete Lagrangian one-forms}. Locally,
\begin{align*}
  & D_{\DEL} L_d (t_{k-1},q_{k-1},t_{k},q_{k},t_{k+1},q_{k+1}) = [ D_{4} L_d
  (t_{k-1},q_{k-1},t_{k},q_{k}) + \\
  & \! D_{2} L_d (t_{k},q_{k},t_{k+1},q_{k+1}) ] dq_k \! + \! [ D_{3} L_d
  (t_{k-1},q_{k-1},t_{k},q_{k}) \! + \! D_{1} L_d
  (t_{k},q_{k},t_{k+1},q_{k+1}) ] dt_k ,
  \\
  & \Theta_{L_d}^+ \hspace*{-1pt} (t_{k},q_{k},t_{k+1},q_{k+1}) \! = \! D_{4} L_d
  (t_{k},q_{k},t_{k+1},q_{k+1}) dq_{k+1} \! + \! D_{3} L_d
  (t_{k},q_{k},t_{k+1},q_{k+1}) dt_{k+1} ,
  \\
  & \Theta_{L_d}^- \hspace*{-1pt} (t_{k},q_{k},t_{k+1},q_{k+1}) \! = \! - D_{2} L_d
  (t_{k},q_{k},t_{k+1},q_{k+1}) dq_{k} - D_{1} L_d
  (t_{k},q_{k},t_{k+1},q_{k+1}) dt_{k} ,
\end{align*}
where  $D_i$ denotes  the differential  with  respect to  the $i$th  variable,
$i=1,\dots,4$.

To ease the exposition, along the paper we consider smooth discrete
Lagrangians.

\noindent {\bf Discrete Hamilton principle}. The discrete variational
principle states that, given fixed end points $(t_0,q_0)$, $(t_N,q_N)$, the
evolution equations extremize $S$.

Otherwise said, we seek discrete paths $c \in C_d$ which are critical
points of the discrete action, $dS (c) \cdot \delta c = 0$ for all
variations $\delta c \in T_c C_d$ with $\delta c (0) = 0 = \delta c
(N)$.  From Theorem~\ref{th:discrete-action}, we get the extended
discrete Euler-Lagrange (EDEL) equations,
\begin{align}
  D_{\DEL} L_d (t_{k-1},q_{k-1},t_{k},q_{k},t_{k+1},q_{k+1}) = 0 \, ,
  \quad 1 \le k \le N-1 \, ,
\end{align}
which can be equivalently written us
\begin{align}
  D_2 L_d(t_k,q_k,t_{k+1},q_{k+1}) + D_4
  L_d(t_{k-1},q_{k-1},t_{k},q_{k}) = 0 \, , \label{eq:EDEL-I}
  \\
  D_1 L_d(t_k,q_k,t_{k+1},q_{k+1}) + D_3L_d(t_{k-1},q_{k-1},t_k,q_{k})
  = 0 \, . \label{eq:EDEL-II}
\end{align}

If we define the discrete energies of the system to be
\begin{align*}
  E_{L_d}^+ (t_k,q_k,t_{k+1},q_{k+1}) & = - D_3
  L_d(t_k,q_k,t_{k+1},q_{k+1}) \, ,
  \\
  E_{L_d}^- (t_k,q_k,t_{k+1},q_{k+1}) &= D_1
  L_d(t_{k},q_{k},t_{k+1},q_{k+1}) \, ,
\end{align*}
then equation~\eqref{eq:EDEL-II} can be simply written as
\[
E_{L_d}^+ (t_{k-1},q_{k-1},t_{k},q_{k}) = E_{L_d}^-
(t_k,q_k,t_{k+1},q_{k+1}) \, ,
\]
which reflects the evolution of the discrete energies.

Under appropriate regularity conditions on the discrete Lagrangian $L_d$
(see~\cite{MaWe}), the DEL equations induce an \emph{extended discrete
  Lagrangian map} $\Phi: \ov{Q} \times \ov{Q} \rightarrow \ov{Q} \times
\ov{Q}$, $(t_{k-1},q_{k-1},t_{k},q_{k}) \mapsto (t_k,q_k,t_{k+1},q_{k+1})$.
The basic geometric properties concerning extended variational integrators
derived from the EDEL equations are the following,

{\bf Symplecticity}: consider the  restricted discrete action $\hat{S}: \ov{Q}
\times \ov{Q} \rightarrow \real$,
\[
\hat{S} (t_0,q_0,t_1,q_1) = S (c) \, ,
\]
where $c \in C_d$ is the unique solution of the EDEL equations
satisfying $c(0) = (t_0,q_0)$, $c(1) = (t_1,q_1)$. From
Theorem~\ref{th:discrete-action}, we compute $d\hat{S} =
(\Phi^{N-1})^* \Theta_{L_d}^+ - \Theta_{L_d}^-$, and then
\[
(\Phi^{N-1})^* \Omega_{L_d} = \Omega_{L_d} \, ,
\]
where $\Omega_{L_d}$ is the extended discrete Lagrangian two-form,
$\Omega_{L_d} = - d \Theta_{L_d}^+ = - d \Theta_{L_d}^-$. Therefore,
extended variational integrators are
symplectic~\cite{KaMaOr,LeMa,MaWe}.

{\bf Extended Noether's theorem}: Let $G$ be a Lie group acting on
$\ov{Q}$, $\psi: G \times \ov{Q} \rightarrow \ov{Q}$, and consider its
diagonal extension to $\ov{Q} \times \ov{Q}$,
\[
\setlength{\arraycolsep}{2pt}
\begin{array}{rccl}
  \Psi: & G \times \ov{Q} \times \ov{Q} & \longrightarrow & \ov{Q} \times
  \ov{Q} 
  \\
  & (g,t_0,q_0,t_1,q_1) & \longmapsto & (\psi (g,t_0,q_0), \psi (g,t_1,q_1))
  \, .  
\end{array} 
\]
The discrete Lagrangian $L_d$ is \emph{$G$-invariant} if $L_d (\Psi
(g,t_0,q_0,t_1,q_1)) = L_d (t_0,q_0,t_1,q_1)$, for all $g \in G$, $(t_0,q_0),
(t_1,q_1) \in \ov{Q}$.  The discrete Lagrangian $L_d$ is \emph{infinitesimally
  invariant} if $\langle dL_d, \xi_{\ov{Q} \times \ov{Q}} \rangle = 0$,
$\forall \xi \in \g$, where $\xi_{\ov{Q} \times \ov{Q}} (t_0,q_0,t_1,q_1) =
(\xi_{\ov{Q}}(t_0,q_0),\xi_{\ov{Q}}(t_1,q_1))$ denotes the infinitesimal
generator of $\Psi$ associated with $\xi$.  Clearly, invariant Lagrangians are
also infinitesimally invariant.  Using $dL_d = \Theta_{L_d}^+ -
\Theta_{L_d}^-$, one sees that an infinitesimally invariant Lagrangian defines
a canonical \emph{discrete momentum map},
\[
\setlength{\arraycolsep}{2pt}
\begin{array}{rccl}
  J_{L_d}: & \ov{Q} \times \ov{Q} & \longrightarrow & \g^* \\
  & (t_0,q_0,t_1,q_1) & \longmapsto & J_{L_d} (t_0,q_0,t_1,q_1):
  \begin{array}[t]{rcl}
    \g & \rightarrow & \real \\
    \xi & \mapsto & \Theta_{L_d}^+ \cdot \xi_{\ov{Q} \times \ov{Q}} =
    \Theta_{L_d}^- \cdot \xi_{\ov{Q} \times \ov{Q}} \, .  
  \end{array}
\end{array}
\]
If $L_d$ is $G$-invariant, then it can be easily seen that $\Psi_g^*
\Theta_{L_d}^{\pm} = \Theta_{L_d}^{\pm}$.  This implies that $J_d$ is
$\Ad$-equivariant.  A second fundamental fact is that extended variational
integrators preserve momentum~\cite{KaMaOr,LeMa,MaWe}, i.e.  $J_{L_d} \circ
\Phi = J_{L_d}$.

{\bf Energy conservation for autonomous discrete Lagrangians}: a
discrete Lagrangian is called \emph{autonomous} if it is invariant
with respect to the additive action of $\real$ on the time component
of $\ov{Q}$, $\psi : \real \times \ov{Q} \rightarrow \ov{Q}$, $\psi
(s,(t,q))=(s+t,q)$. The associated discrete momentum map is given by
$J_{L_d} (t_0,q_0,t_1,q_1) = - E_{L_d}^+ (t_0,q_0,t_1,q_1) dt_1 =
-E_{L_d}^- (t_0,q_0,t_1,q_1) dt_0$. Noether's theorem thus gives
\[
E_{L_d}^+ (t_k,q_k,t_{k+1},q_{k+1}) = E_{L_d}^+
(t_{k-1},q_{k-1},t_{k},q_{k}) \, ,
\]
or equivalently, $E_{L_d}^- (t_k,q_k,t_{k+1},q_{k+1}) = E_{L_d}^-
(t_{k-1},q_{k-1},t_{k},q_{k})$, i.e. the discrete energy is conserved
by the extended variational integrators derived from an autonomous
Lagrangian~\cite{KaMaOr,LeMa,MaWe}.

\section{A discrete Lagrange-d'Alembert principle for nonautonomous
  systems}\label{se:EDLA}

In this section, we propose a discrete version of the Lagrange-d'Alembert
principle for nonautonomous discrete systems.  We start by defining what we
understand by an extended discrete nonholonomic Lagrangian system,

\begin{definition}
  An \emph{extended discrete nonholonomic Lagrangian system} is a
  triple $(L_d,\D_d,\D)$, where $L_d : \ov{Q} \times \ov{Q}
  \rightarrow \real$ is the discrete Lagrangian, ${\cal D}_d \subset
  \ov{Q} \times \ov{Q}$ is the discrete constraint space and ${\cal
    D}$ is the constraint distribution on $\ov{Q}$.  In addition,
  $\D_d$ has the same dimension as $\D$ and is such that $(t,q,t,q)
  \in \D_d$ for all $(t,q) \in \ov{Q}$.
\end{definition}

Notice that the unconstrained discrete mechanics (cf.
Section~\ref{se:variational-integrators}) can also be seen within this
framework, where ${\cal D}=T \ov{Q}$ and ${\cal D}_d = \ov{Q} \times
\ov{Q}$.

\begin{remark}
  The motivation for this notion of extended discrete nonholonomic
  Lagrangian system is the following.  When dealing with unconstrained
  systems, given fixed end points $(t_0,q_0)$, $(t_N,q_N)$, one
  extremizes the action sum $S$ with respect to all possible discrete
  paths. This means that at each point $(t,q) \in \ov{Q}$, the allowed
  variations are the whole tangent space $T_{(t,q)} \ov{Q}$.  However,
  in the nonholonomic case, one must restrict the allowed variations
  at each point: these will be exactly given by the distribution of
  feasible velocities $\D$.  On the other hand, the discrete
  constraint space $\D_d$ will impose certain constraints on the
  solution sequence $\{ (t_k,q_k) \}$.
\end{remark}

Here, we will only consider constraints which do not impose conditions
on the time velocities, i.e.  $\tau_* (\D) = T \real$, where $\tau:
\ov{Q} = \real \times Q \rightarrow \real$ is the projection onto the
first factor, although most of the discussion can be also carried out
in broader terms.  The \emph{constrained discrete path space} is the
set of extended discrete paths which verify the discrete constraints,
\[
\tilde{C}_d = \{ c \in C_d \; | \; c (k) \in \D_d \, , \; 0 \le k \le N \} \, ,
\]
and the set of allowed variations is given by
\[
\V_{d} = \{ \delta c \in TC_d \; | \; \delta c (k) \in \D_{c(k)} \, ,
\; 0 \le k \le N \} \, .
\]

\noindent {\bf Discrete Lagrange-d'Alembert principle}. Given fixed end points
$(t_0,q_0)$ and $(t_N,q_N)$, the discrete Lagrange-d'Alembert principle
consists of extremizing the extended action sum $S$ among the variations in
$\V_d$ and such that the solution sequence belongs to $\tilde{C}_d$.

Otherwise said, we seek discrete paths $c \in \tilde{C}_d$ such that
$dS (c) \cdot \delta c = 0$, for all $\delta c \in \V_d$, with $\delta
c (0) = 0 = \delta c (N)$. Using Theorem~\ref{th:discrete-action}, we
get
\begin{align*}
  0 = dS (c) \cdot \delta c = \sum_{k=1}^{N-1} D_{\DEL} L_d
  (t_{k-1},q_{k-1},t_{k},q_{k},t_{k+1},q_{k+1}) \cdot (\delta t_{k},
  \delta q_{k}) \, ,
\end{align*}
for all $(\delta t_{k}, \delta q_{k}) \in \D_{(t_k,q_k)}$, $1 \le k
\le N-1$. Hence, the extended discrete Lagrange-d'Alembert (EDLA)
equations read
\begin{equation*}
  \left\{
    \begin{array}{l}
      D_{\DEL} L_d (t_{k-1},q_{k-1},t_{k},q_{k},t_{k+1},q_{k+1}) \in
      \D_{(t_k,q_k)}^o \, , 
      \\
      (t_k,q_k,t_{k+1},q_{k+1}) \in \D_d \, , 
    \end{array}
  \right.
\quad 1 \le k \le N-1 \, ,
\end{equation*}
where $\D^o$ denotes the annihilator of $\D$. Let $\omega_d^a: \ov{Q} \times
\ov{Q} \rightarrow \real$, $a \in \{ 1,\dots, m \}$, be smooth functions whose
annihilation defines locally $\D_d$, and let $\omega^a : \ov{Q} \rightarrow
T^*\ov{Q}$, $a \in \{ 1,\dots, m \}$ be one-forms on $\ov{Q}$ locally defining
$\D^o$.  Since $\tau_* (\D) = T \real$, the latter ones are of the form
$\omega (t,q) = (0, \omega (t,q))$, where with a slight abuse of notation we
denote in the same way the component of the one-form in $T^*Q$ and the
one-form itself. The EDLA equations can then be written as,
\begin{equation}\label{eq:EDLA}
  \left\{
    \setlength{\arraycolsep}{2pt}
    \begin{array}{l}
      D_1 L_d(t_k,q_k,t_{k+1},q_{k+1}) + D_3 L_d(t_{k-1},q_{k-1},t_{k},q_{k})
      = 0 \, , \\ 
      D_2 L_d(t_k,q_k,t_{k+1},q_{k+1}) + D_4 L_d(t_{k-1},q_{k-1},t_k,q_{k}) =
      \lambda_a \omega^a (t_k,q_k) \, , \\ 
      \omega_d^a (t_k,q_k,t_{k+1},q_{k+1}) = 0 \, .
    \end{array}
  \right.
\end{equation}
Notice that the discrete Lagrange-d'Alembert principle is not truly
variational, in the sense that it does not correspond to the
extremization of any action sum.  This is in accordance with the
nature of its continuous counterpart.  Alternatively, we will refer to
the EDLA algorithm~\eqref{eq:EDLA} as a \emph{nonholonomic
  integrator}, by analogy with the unconstrained case.

\begin{remark}[Well-posedness of the discrete problem]
  As it is also the case in unconstrained discrete
  mechanics~\cite{KaMaOr}, the existence of solutions for the extended
  equations is not always guaranteed. If the mapping
  \[
  \setlength{\arraycolsep}{2pt}
  \begin{array}{rcl}
    \D_d \times \real^m & \rightarrow & T^*\ov{Q} \\
    (t_0,q_0,t_1,q_1,\lambda) & \mapsto & (t_0,q_0,D_1
    L_d(t_0,q_0,t_1,q_1),-D_2 
    L_d(t_0,q_0,t_1,q_1) + \lambda_a \omega^ a(t_0,q_0)) \, ,
  \end{array}
  \]
  is a local diffeomorphism, then for a pair $(t_{k-1},q_{k-1})$,
  $(t_k,q_k)$, there exists $(t_{k+1},q_{k+1})$ verifying the EDLA
  equations~\eqref{eq:EDLA}.  The problem now arises from the fact
  that $t_{k+1}>t_k$ is not guaranteed, and therefore one might obtain
  inconsistent solutions.  Nevertheless, one can ensure that, for
  specific choices of discrete Lagrangians~\cite{KaMaOr} of natural
  (kinetic minus potential energy) systems, this situation does not
  occur away from points where the discrete energy is near zero.
\end{remark}

In the remainder of the paper, we assume that the EDLA
equations~\eqref{eq:EDLA} are well-posed and therefore induce an
\emph{extended discrete Lagrange-d'Alembert map} $\Phi: \ov{Q} \times \ov{Q}
\rightarrow \ov{Q} \times \ov{Q}$, $(t_{k-1},q_{k-1},t_{k},q_{k}) \mapsto
(t_k,q_k,t_{k+1},q_{k+1})$. The actual implementation of the EDLA algorithm
can be carried out building on the discussion in~\cite{CoMa,KaMaOr,MaWe}.

\section{Geometric properties}\label{se:EDLA-geometric-properties}

In this section, we examine the geometric properties of the
integrators derived from the discrete Lagrange-d'Alembert principle
proposed above. It is important to keep in mind that the continuous
flow of a nonholonomic Lagrangian problem does not have the same
properties as the unconstrained flow~\cite{Co}: on the one hand, the
Poincar\'e-Cartan form $\Omega_L$ is no longer preserved in general.
On the other hand, the action of a symmetry Lie group does not
generally give rise to momentum conserved quantities.  However, the
nonholonomic flow does enjoy some nice geometric properties with
respect to these objects, which will guide our study of the
corresponding discrete mechanics.

{\bf Symplectic form}: Consider the restricted action $\tilde{S}:
\ov{Q} \times \ov{Q} \rightarrow \real$,
\[
\tilde{S} (t_0,q_0,t_1,q_1) = S(c) \, ,
\]
where $c$ is the unique solution of the EDLA equations satisfying
$c(0) = (t_0,q_0)$, $c(1) = (t_1,q_1)$. Using
Theorem~\ref{th:discrete-action} with $N=2$, we compute
\[
d \tilde{S} = \lambda_a \omega^a (t_1,q_1) + \Phi^* \Theta_{L_d}^+ -
\Theta_{L_d}^- \, ,
\]
and therefore conclude that $\Phi^* \Omega_{L_d} = \Omega_{L_d} + d
\beta_d$, with $\beta_d \in \D^o$.  Note that this is the discrete
version of the behavior of the nonautonomous continuous flow with
respect to the Poincar\'e-Cartan two-form, $\Lc_{X} \Omega_L = d
\beta$, with $\beta \in (\D^v)^o$ (see~\cite{Co}).

{\bf Nonholonomic momentum map}: Assume that the extended discrete
nonholonomic Lagrangian system $(L_d,\D_d,\D)$ is invariant under the
(diagonal) action of a Lie group $G$ on $\ov{Q}$, that is, all the
three elements are $G$-invariant. Let $\V$ denote the bundle of
vertical vectors with respect to the canonical projection $\pi :
\ov{Q} \rightarrow \ov{Q}/G$,
\[
\V_{(t,q)} = \{ \xi_{\ov{Q}} (t,q) \; | \; \xi \in \g \} \, .
\]
Among these symmetry directions, we are interested in selecting those
ones which are also compatible with the nonholonomic constraints, that
is,
\[
\g_{(t,q)} = \{ \xi \in \g \; | \; \xi_{\ov{Q}}(t,q) \in \D_{(t,q)} \}
\, .
\]
Let $\g^{\D}$ denote the (generalized) vector bundle over $\ov{Q}$ whose fiber
at $(t,q)$ is  given by $\g_{(t,q)}$. We now  define the discrete nonholonomic
momentum map as,
\[
\setlength{\arraycolsep}{2pt}
\begin{array}{rccl}
  J_{d}^{nh} : & \ov{Q} \times \ov{Q} & \longrightarrow & (\g^\D)^* \\
  & (t_0,q_0,t_1,q_1) & \longmapsto & J_{d}^{nh} (t_0,q_0,t_1,q_1):
  \begin{array}[t]{rcl}
    \g^\D & \rightarrow & \real \\
    \xi & \mapsto & \langle J_d (t_0,q_0,t_1,q_1) , \xi \rangle \, .
  \end{array}
\end{array}
\]
Note that this mapping is just the restriction of the usual discrete
momentum map to the fiber bundle $\g^\D$.  Now, take a
$C^{\infty}$-section of $\g^\D \rightarrow \ov{Q}$, that is, a mapping
$\tilde{\xi}$ which for each $(t,q) \in \ov{Q}$ gives us a symmetry
direction $\tilde{\xi} (t,q)$ whose associated fundamental vector
field lies in the constraint distribution.

\begin{proposition}\label{prop:discrete-momentum}
  Assume that $(L_d,\D_d,\D)$ is invariant under the action of $G$.
  Then, the discrete time evolution of the nonholonomic momentum map
  is governed by the \emph{discrete momentum equation},
  \begin{align}\label{eq:discrete-momentum-equation}
    & \langle J_d^{nh} (t_1,q_1,t_2,q_2), \tilde{\xi} \rangle -
    \langle J_d^{nh}
    (t_0,q_0,t_1,q_1), \tilde{\xi} \rangle \nonumber \\
    & \quad \qquad = \langle \Theta_{L_d}^+ (t_1,q_1,t_2,q_2) ,
    (\tilde{\xi}(t_2,q_2) - \tilde{\xi} (t_1,q_1))_Q (t_2,q_2) \rangle
    \, .
  \end{align}
\end{proposition}

\begin{proof}
  The Lie group $G$ acts on $C_d$ by means of the pointwise action.
  Then,
  \[
  \langle d S (c) , \xi_{C_d} (c) \rangle = \sum_{k=0}^{N-1} \langle
  dL_d , \xi_{\ov{Q} \times \ov{Q}} \rangle = 0 \, .
  \]
  On the other hand, since the space $\D_d$ is $G$-invariant,
  $\tilde{C}_d$ is preserved by the group action. All this, together
  with the invariance of $\D$, implies that the solutions to the EDLA
  equations~\eqref{eq:EDLA} are preserved by $G$, i.e., $\Phi \circ
  \Psi_g = \Psi_g \circ \Phi$.
  
  Let $(t_0,q_0,t_1,q_1) \in \ov{Q} \times \ov{Q}$ and consider the
  corresponding solution to the EDLA equations. Take $N=2$ and then,
  \begin{multline*}
    0 = \langle dS (c) , \xi_{C_d} (c) \rangle = \langle d \tilde{S}
    (t_0,q_0,t_1,q_1) , \xi_{\ov{Q} \times \ov{Q}} (t_0,q_0,t_1,q_1) \rangle \\
    = \left\langle \lambda_a \omega^a (t_1,q_1) + \Phi^*
      \Theta_{L_d}^+ - \Theta_{L_d}^- , \xi_{\ov{Q} \times \ov{Q}}
      (t_0,q_0,t_1,q_1) \right\rangle\, .
  \end{multline*}
  Now, if $\xi_{\ov{Q}} (t_1,q_1)$ belongs to $\D_{(t_1,q_1)}$, we
  deduce that
  \[
  \langle \Theta_{L_d}^+(t_1,q_1,t_2,q_2) , \xi_{\ov{Q}} (t_2,q_2)
  \rangle = \langle \Theta_{L_d}^- (t_0,q_0,t_1,q_1) , \xi_{\ov{Q}}
  (t_0,q_0) \rangle \, .
  \]
  Finally,
  \begin{multline*}
    \langle J_d^{nh} (t_1,q_1,t_2,q_2), \tilde{\xi} \rangle - \langle
    J_d^{nh} (t_0,q_0,t_1,q_1), \tilde{\xi} \rangle
    \\
    = \langle \Theta_{L_d}^+ (t_1,q_1,t_2,q_2) ,
    (\tilde{\xi}(t_2,q_2))_{\ov{Q}} (t_2,q_2) \rangle - \langle
    \Theta_{L_d}^- (t_0,q_0,t_1,q_1) ,
    (\tilde{\xi}(t_1,q_1))_{\ov{Q}} (t_0,q_0) \rangle \\
    = \langle \Theta_{L_d}^+ (t_1,q_1,t_2,q_2) ,
    (\tilde{\xi}(t_2,q_2))_{\ov{Q}} (t_2,q_2) \rangle - \langle
    \Theta_{L_d}^+ (t_1,q_1,t_2,q_2) , (\tilde{\xi}(t_1,q_1))_{\ov{Q}}
    (t_2,q_2) \rangle \, ,
  \end{multline*}
  which is the desired result.
\end{proof}

A distinguished class of sections of the bundle $\g^\D$ is formed by
the constant ones, $\tilde{\xi} (q) = \xi$. They correspond to
elements $\xi$ of the Lie algebra which always are compatible with the
constraints, that is, $\xi_{\ov{Q}} (t,q) \in \D_{(t,q)}$, for all
$(t,q) \in \ov{Q}$. These special elements are called \emph{horizontal
  symmetries} in the literature of nonholonomic
mechanics~\cite{Bl,BlKrMaMu,Co}.

\begin{corollary}\label{cor:horizontal-symmetry}
  If $\xi \in \g$ is a horizontal symmetry, then the associated
  component of the discrete nonholonomic momentum is preserved by the
  EDLA algorithm.
\end{corollary}

\begin{proof}
  It is immediate from~\eqref{eq:discrete-momentum-equation}, since in
  this case $\tilde{\xi}(t_2,q_2) - \tilde{\xi} (t_1,q_1) = 0$, and
  hence $\langle J_d^{nh} (t_1,q_1,t_2,q_2), \tilde{\xi} \rangle =
  \langle J_d^{nh} (t_0,q_0,t_1,q_1), \tilde{\xi} \rangle$.
\end{proof}

{\bf Nonholonomic Chaplygin systems}: It may also happen that the
generalized bundle $\g^\D$ over $Q$ is trivial, that is,
$\g_{(t,q)}=0$ for all $(t,q) \in \ov{Q}$. In this case, there is no
nonholonomic momentum map and hence we must look for different
geometric properties of the flow other than
Proposition~\ref{prop:discrete-momentum}.  Under the additional
hypothesis $\D + \V = T\ov{Q}$ (\emph{dimensional assumption}, cf.
\cite{Bl,BlKrMaMu}), we deduce that $\D$ complements $\V$ in the
tangent bundle of $\ov{Q}$, and therefore constitutes the horizontal
space of a principal connection.  We denote its associated connection
one-form by $\A :T\ov{Q} \rightarrow \g$.  This class of nonholonomic
systems are called generalized Chaplygin systems~\cite{CaCoLeMa}.
It is known that, after the reduction by the action of the Lie group,
these systems give rise to an unconstrained system subject to an
external force of gyroscopic type.  In the following, we show that the
discrete mechanics also shares this feature.

Assume that the discrete constraint space $\D_d$ and the action are
such that $T\D_d + \V \times \V = T\ov{Q}$ (an hypothesis that we term
\emph{discrete dimensional assumption}).  Let $\pi: \ov{Q} \rightarrow
\ov{Q}/G$ be the canonical projection, and consider the map
\[
\setlength{\arraycolsep}{2pt}
\begin{array}{rccl}
  \nu: & \D_d / G & \longrightarrow & \ov{Q} /G \times \ov{Q} /G \\
  & [(t_0,q_0,t_1,q_1)] & \longmapsto & (\pi (t_0,q_0), \pi (t_1,q_1)) \, .
\end{array}
\]
Note that both spaces have the same dimension due to the definition of
$\D_d$ and the dimensional assumption. Indeed, $\dim \D_d = \dim \D =
2 \dim \ov{Q} - \dim G$. On the other hand, if $\rho:\D_d \rightarrow
\D_d/G$ denotes the projection from $\D_d$ to its reduced space, then
one can verify that $\ker \rho_* \subset T \D_d \cap (\V \times \V)$.
By a dimensional argument, we conclude that $\ker \rho_* = T \D_d \cap
(\V \times \V)$, and therefore $\nu$ is a local diffeomorphism.

We say that $\D_d$ is \emph{right-rigid with respect to the $G$-action
  $\psi$} if the following property holds: given $(t_0,q_0,t_1,q_1)
\in \D_d$ and $g \in G$, if $(t_0,q_0,g(t_1,q_1)) \in \D_d$, then
$g=e$ (where we are using the abbreviated notation $g(t_1,q_1) = \psi
(g,t_1,q_1))$).  Clearly, if $\D_d$ is right-rigid and invariant under
the diagonal action, it is also left-rigid.  Intuitively, the notion
of right-rigidity (resp. left-rigidity) means that $\D_d$ is not
invariant under the action $\Id \times \psi : G \times \ov{Q} \times
\ov{Q} \rightarrow \ov{Q} \times \ov{Q}$, $(\Id \times \psi)
(g,t_0,q_0,t_1,q_1) = (t_0,q_0,\psi(g,t_1,q_1))$ (resp.  $\psi \times
\Id$).

\begin{proposition}\label{prop:global}
  Let $(L_d,\D_d,\D)$ be $G$-invariant.  Assume that the generalized
  bundle $\g^{\D}$ on $\ov{Q}$ is trivial and that the discrete
  dimensional assumption holds.  Then, if $\D_d$ is right-rigid, the
  local diffeomorphism $\nu$ is global.
\end{proposition}
\begin{proof}
  We use the abbreviated notation $\ov{q}=(t,q) \in \ov{Q}$.  Take
  $[(\ov{q}_0,\ov{q}_1)]$, $[(\tilde{\ov{q}}_0, \tilde{\ov{q}}_1)] \in
  \D_d /G$ such that $\nu [(\ov{q}_0,\ov{q}_1)] = \nu
  [(\tilde{\ov{q}}_0, \tilde{\ov{q}}_1)]$. Then there exist $g_0$,
  $g_1 \in G$ such that $\ov{q}_0 = g_0 \tilde{\ov{q}}_0$, $\ov{q}_1 =
  g_1 \tilde{\ov{q}}_1$.  Since $(\tilde{\ov{q}}_0,\tilde{\ov{q}}_1)
  \in \D_d$, then $(g_0 \tilde{\ov{q}}_0, g_0 \tilde{\ov{q}}_1) \in
  \D_d$ by $G$-invariance.  Alternatively, we have
  $(\ov{q}_0,\ov{q}_1) \in \D_d$ and, at the same time, $(g_0
  \tilde{\ov{q}}_0, g_0 \tilde{\ov{q}}_1) = (\ov{q}_0,g_0 g_1^{-1}
  \ov{q}_1) \in \D_d$.  Now, by rigidity, we conclude $g_0 = g_1$, and
  hence $[(\ov{q}_0,\ov{q}_1)] =
  [(\tilde{\ov{q}}_0,\tilde{\ov{q}}_1)]$.
\end{proof}

Therefore, under the global identification provided by $\nu$, we can
define a reduced discrete Lagrangian $L_d^* : \ov{Q} /G \times \ov{Q}
/G \rightarrow \real$, $L_d^* (\ov{r}_k,\ov{r}_{k+1}) = \ell_d
(\nu^{-1}(\ov{r}_k,\ov{r}_{k+1}))$, where $\ell_d : (\ov{Q} \times
\ov{Q}) / G \rightarrow \real$ is the reduction of $L_d$ to $(\ov{Q}
\times \ov{Q}) / G$, and we regard $\D_d /G$ as a submanifold of
$(\ov{Q} \times \ov{Q}) / G$.  Locally, if we identify $\ov{Q}$ with
$\ov{Q}/G \times G$, $(t,q)=(\ov{r},g)$, then we can take local
coordinates $(\ov{r}_0,\ov{r}_1,f_{0,1}) \in \ov{Q}/G \times \ov{Q}/G
\times G$ on $(\ov{Q}\times \ov{Q})/G$. In this way, the projection
$\ov{Q} \times \ov{Q} \rightarrow (\ov{Q} \times \ov{Q})/G$ reads
$(\ov{r}_0,g_0,\ov{r}_1,g_1) \mapsto
(\ov{r}_0,\ov{r}_1,f_{0,1}=g_0^{-1}g_1)$.  Moreover, when regarding
$\D_d /G$ as contained in $(\ov{Q} \times \ov{Q}) / G$, we have that
$f_{0,1} = f_{0,1} (\ov{r}_0, \ov{r}_1)$ for
$(\ov{r}_0,\ov{r}_1,f_{0,1}) \in \D_d /G$.  Finally, if the $G$-action
acts trivially on the time component of $\ov {Q}$, we can further
write $\ov{Q}/G = \real \times Q/G $, $\ov{r} = (t,r)$.  Now, we are
in a position to state the following result.

\begin{proposition}\label{prop:REDLA}
  Under the hypothesis of Prop.~\ref{prop:global}, assume $\psi$ acts
  trivially on the time component of $\ov {Q}$.  Then, the solutions
  of the EDLA equations project onto the solutions of the reduced
  extended discrete Lagrange-d'Alembert (REDLA) equations,
  \[
  \setlength{\arraycolsep}{2pt}
  \begin{array}{lll}
    D_1 L^*_d(t_k,r_k,t_{k+1},r_{k+1}) + D_3 L^*_d(t_{k-1},r_{k-1},t_{k},r_{k})
    &=& 0 \, , \\ 
    D_2 L^*_d(t_k,r_k,t_{k+1},r_{k+1}) + D_4 L^*_d(t_{k-1},r_{k-1},t_k,r_{k})
    &=& F^-(t_k,r_k,t_{k+1},r_{k+1}) \\ 
    && + F^+(t_{k-1},r_{k-1},t_{k},r_{k}) \, ,
  \end{array}
  \]
  where the expression of the forces in bundle coordinates is given by
  \begin{align*}
    F^-(\ov{r}_k,\ov{r}_{k+1}) &= \frac{\partial \ell_d}{\partial
      f_{k,k+1}} \left( \frac{\partial f_{k,k+1}}{\partial
        r_k}(\ov{r}_k,\ov{r}_{k+1}) -
      R_{f_{k,k+1}(\ov{r}_k,\ov{r}_{k+1})}
      (\A_{\text{loc}}(\ov{r}_k)(\cdot)) \right) \, ,
    \\
    F^+(\ov{r}_{k-1},\ov{r}_{k}) &= \frac{\partial \ell_d}{\partial
      f_{k-1,k}} \left( \frac{\partial f_{k-1,k}}{\partial
        r_k}(\ov{r}_{k-1},\ov{r}_k) +
      L_{f_{k-1,k}(\ov{r}_{k-1},\ov{r}_k)} (\A_{\text{loc}}
      (\ov{r}_k)(\cdot)) \right) \, ,
  \end{align*}
  where $\A_{\text{loc}}(r)$ is the local form of the connection one-form $\A$.
\end{proposition}

This result can be proved using a similar argument to the one carried out
in~\cite{CoMa} for autonomous systems.  The REDLA integrator is an appropriate
version for nonautonomous systems of the \emph{generalized variational
  integrators} developed for systems subject to external forcing
in~\cite{KaMaOrWe}.  This is in accordance with the situation in the
continuous case where, as we mentioned before, the reduction of the Chaplygin
system gives rise to an unconstrained system subject to a gyroscopic external
force.

{\bf Energy conservation for autonomous constrained Lagrangian
  systems}: The discrete system $(L_d,\D_d,\D)$ is \emph{autonomous}
if $L_d$, $\D$ and $\D_d$ are invariant under the additive action of
$\real$ on the time component of $\ov{Q}$. In this case, $\g^\D = \g =
\real$. Consequently, the discrete nonholonomic momentum map coincides
with $J_d$, which, as we have already seen, is given by $J_{L_d}
(t_0,q_0,t_1,q_1) = - E_{L_d}^+ (t_0,q_0,t_1,q_1) dt_1 = -E_{L_d}^-
(t_0,q_0,t_1,q_1) dt_0$. Corollary~\ref{cor:horizontal-symmetry} thus
yields
\[
E_{L_d}^+ (t_k,q_k,t_{k+1},q_{k+1}) = E_{L_d}^+
(t_{k-1},q_{k-1},t_{k},q_{k}) \, .
\]
\begin{proposition}\label{prop:energy}
  If the discrete system $(L_d,\D_d,\D)$ is autonomous, then the EDLA
  algorithm preserves its associated discrete energy.
\end{proposition}

\section{Conclusions}\label{se:conclusions}

We have proposed a discrete version of the Lagrange-d'Alembert
principle for nonautonomous Lagrangian systems with nonholonomic
constraints.  We have studied the geometric properties of the
integrators derived from this principle, paying special attention to
the evolution of the symplectic form and the nonholonomic momentum
map, and the conservation of energy. Future work will be devoted to
develop a numerical error analysis of these integrators making use of
backward error techniques.

\section*{Acknowledgments}

This work was supported by DGICYT grant BFM2001-2272.

\end{document}